\documentclass[art10]{article}
\usepackage{graphicx}
\usepackage{graphpap}
\usepackage[mathscr]{eucal}
\usepackage{amsmath}
\usepackage{amsfonts}
\usepackage{amssymb}
\newcommand{\bc}{\begin{center}}
\newcommand{\ec}{\end{center}}
\newcommand{\be}{\begin{eqnarray}}
\newcommand{\ee}{\end{eqnarray}}

\newcommand{\BOX}{\hspace{0.5cm} $\Box$}
\newcommand{\pr}{{\bf Proof}\hspace{0.3cm}}

\begin{document}

\author{ }
\date{}
%\maketitle
\newtheorem{lemma}{Lemma}
\newtheorem{theorem}{Theorem}
\newtheorem{proposition}{Proposition}
\newtheorem{remark}{Remark}
\newtheorem{corollary}{Corollary}
\newtheorem{open problem}{Open problem}
\newtheorem{con}{Conjecture}
%\Large \textbf{Hasil Penelitian}\normalsize
\bc {\large \textbf{Super edge-magic deficiency of join-product graphs}}\\
\bigskip
\vspace{1cm} A.A.G. Ngurah\footnote{The author was funded by ``Hibah Desentralisasi - Fundamental 2013",  088/SP2H/PDSTRL/K7/KL/III/2013, from the Directorate General of Higher Education, Indonesia.}\\
{\small Department of Civil Engineering\\ Universitas Merdeka Malang\\ Jalan Taman Agung No. 1 Malang, Indonesia\\email: ngurahram67@yahoo.com}\\ Rinovia Simanjuntak\\  {\small Combinatorial Mathematics Research Group\\
Faculty of Mathematics and Natural Sciences\\
Institut Teknologi Bandung\\
Jl. Ganesa 10 Bandung 40132 Indonesia\\email: rino@math.itb.ac.id}
\ec
\begin{quote}
\bigskip
\end{quote}
\begin{center}
\textbf{Abstract}
\end{center}
\begin{quote}
{\small A graph $G$ is called \textit{super edge-magic} if there exists a bijective function $f$
from $V(G) \cup E(G)$ to $\{1, 2, \ldots, |V(G) \cup E(G)|\}$ such that $f(V(G)) = \{1, 2, \ldots, |V(G)|\}$ and $f(x)
+ f(xy) + f(y)$ is a constant $k$ for every edge $xy$ of $G$.
Furthermore,  the \textit{super edge-magic deficiency}
of a graph $G$ is either the minimum nonnegative integer $n$ such that $G \cup nK_1$ is super
edge-magic or $+\infty$ if there exists no such integer.

\emph{Join product} of two graphs is their graph union with additional edges that connect  all vertices of the first graph to each  vertex of the second graph. In this paper, we study the super edge-magic deficiencies of a wheel minus an edge and join products of a path, a star, and a cycle, respectively, with isolated vertices. In general, we show that the join product of a super edge-magic graph with isolated vertices has finite super edge-magic deficiency.}
\end{quote}


\section{Introduction}
All graphs that we consider in this paper are finite and simple. For
most graph theory notions, we refer the reader to Chartrand and
Lesniak's \cite{C}. However, to make this paper reasonably
self-contained, we mention that for a graph $G$, we denote the
vertex and edge sets of graph $G$ by $V(G)$ and $E(G)$,
respectively, and $p = |V(G)|$ and $q = |E(G)|$.

An \textit{edge-magic  labeling} of a graph $G$ is a bijective function $f$
from $V(G) \cup E(G)$ to $\{1, 2, \ldots, p +q\}$ such that  $f(x)
+ f(xy) + f(y)$ is a constant $k$, called a \textit{magic
constant} of $f$, for any edge $xy$ of $G$. An edge-magic labeling
$f$ is called a \textit{super edge-magic labeling} if $f(V(G)) = \{1,
2, \ldots, p\}$. A graph $G$  is called \textit{edge-magic}
(\textit{super edge-magic}) if there exists an edge-magic (super
edge-magic, respectively)  labeling of $G$. The concept of
edge-magic  labeling was first introduced by Kotzig and
Rosa \cite{K} and the super edge-magic labeling was introduced by Enomoto,
Llad\'{o}, Nakamigawa and Ringel \cite{E}. %Recently several papers
%of (super) edge-magic total labeling have been published by several
%authors, for instance see \cite{F3, F4, F5, I, S}.
We mention that an equivalent concept to the one of super edge-magic graphs had already
appeared in the literature under the name of strongly indexable graphs \cite{Ac}. Although the definitions
of super edge-magic graphs and strongly indexable graphs were introduced from different points of
view, they turn out to be equivalent.

In \cite{K},  Kotzig and Rosa proved that for every graph $G$ there
exists an edge-magic graph $H$ such that $H \cong G \cup nK_1$ for
some nonnegative integer $n$. This fact motivated them to define the
concept of edge-magic deficiency of a graph. The \textit{edge-magic
deficiency} of a graph $G$, $\mu(G)$, is defined as the minimum nonnegative
integer $n$ such that $G \cup nK_1$ is  edge-magic. They also proved that every graph has finite edge-magic deficiency.
Motivated by Kotzig and Rosa's concept, Figueroa-Centeno \textit{et al.} \cite{F1}
defined a similar concept for super edge-magic labelings. The
\textit{super edge-magic deficiency} of a graph $G$, $\mu_s(G)$, is either
the minimum nonnegative integer $n$ such that $G \cup nK_1$ is
super edge-magic  or $+\infty$ if there exists no such integer.
As a direct consequence of the above two definitions, the inequality
 $\mu(G) \leq \mu_s(G)$ holds for every graph $G$.
%Unlike the edge-magic deficiency, not all graphs have finite super
%edge-magic deficiency. Examples of such graphs can be found in
%\cite{F1}.

Some authors have studied the super
edge-magic deficiency of some classes of graphs. Figueroa-Centeno \textit{et al.} in two separate papers \cite{F1,
F2} investigated super edge-magic deficiencies of complete graphs, complete bipartite graphs $K_{2,m}$, some classes of
forests with two components, 1-regular graphs, and 2-regular graphs. Ngurah \textit{et al.} \cite{N1, N2} studied the super edge-magic deficiency of
some classes of chain graphs, wheels, fans, double fans, and disjoint union of particular type of complete bipartite graphs. Recently, Ahmad and Muntaner-Battle \cite{A} studied the super edge-magic deficiency of several classes of unicyclic graphs.
The authors refer the reader to the survey paper by Gallian \cite{Ga} for some of the
latest developments in these and other types of graph labelings.

In this paper, we study the super edge-magic deficiencies of a wheel minus an edge and join products of a path, a star, and a cycle, respectively, with isolated vertices. In proving the main results, the following two lemmas  will be used frequently. The first
lemma  characterizes super edge-magic graphs and the second
gives necessary conditions for the existence of super edge-magic graphs.

\begin{lemma} \cite{F}\label{l1}
A graph $G$ with $p$ vertices and $q$ edges is super edge-magic if
and only if there exists a bijective function $f :V(G) \rightarrow
\{1, 2, \ldots, p\}$ such that the set $S = \{f(x) + f(y) : xy \in
E(G)\}$ consists of $q$ consecutive integers. In such a case, $f$
extends to a super edge-magic total labeling of $G$ with the magic
constant $k = p + q + s$, where $s = min(S)$.
\end{lemma}

\begin{lemma} \cite{E}\label{l2}
If a graph $G$ with $p$ vertices and $q$ edges is super
edge-magic, then $q \leq 2p - 3$.
\end{lemma}

%------------------------------------------------------------------------
\section{Super edge-magic deficiency of a wheel minus an edge}

In this section, we consider the super edge-magic deficiency of $W_{n} \cong C_n +
K_1$, $n \geq 3$, minus an edge. We shall denote vertex-set of $W_n$, $V(W_n) = \{c\} \cup \{x_1, x_2, x_3, \ldots, x_n\}$, and edge-set $E(W_n) = \{cx_i : 1 \leq i \leq n\} \cup \{x_ix_{i+1} : 1 \leq i \leq n-1\} \cup \{x_nx_1\}$. We shall call an edge $x_ix_{i+1}$ as a \emph{rim} and an edge $cx_i$ as a \emph{spoke}.
Let us consider the graph $H_n \cong W_{n} - \{e\}$ with order $n+1$ and size $2n-1$. It is interesting to mention that $H_n \cong W_{n} - \{e\}$ is a graph attaining $|E(H_n)| = 2|V(H_n)|- 3$, which is the upper bound of condition in Lemma \ref{l2}.
If the edge $e$ is a rim of $W_{n}$, then $H_n$ is a fan $F_n$ whose super edge-magic deficiency has been studied by Ngurah \textit{et al.} \cite{N1}. They determined the super edge-magic deficiency of $F_n$ for small values of
$n$ and provided upper and lower bounds for general $n$. Here, we consider the super edge-magic
deficiency of $H_n \cong W_{n} - \{e\}$, where $e$ is a spoke of $W_n$. We shall use the following notations for vertex and edge sets: $V(H_n) = \{c\} \cup \{x_1, x_2, x_3, \ldots, x_n\}$, and edge-set $E(H_n) = \{cx_i : 2 \leq i \leq n\} \cup \{x_ix_{i+1} : 1 \leq i \leq n-1\} \cup \{x_nx_1\}$.

Our first result gives the only two super edge-magic labelings for $H_n$.

\begin{theorem} \label{T3}
Let $n \geq 3$ be an integer. The graph $H_n \cong W_{n} - \{e\}$
is super edge-magic if and only if $n \leq 4$.
\end{theorem}
\pr First, we show that $H_n$ is super edge-magic for $n \leq 4$.
Label the vertices $(c; x_1, x_2, x_3)$ and $(c; x_1, x_2, x_3,
x_4)$ with $(1; 4, 3, 2)$ and $(2; 3, 1, 4, 5)$, respectively. This
vertex labeling extends to a super edge-magic labeling of $H_3$ and
$H_4$, respectively.

For the necessity, assume that $H_n$ is super edge-magic with a super
edge-magic labeling $f$ for every integer $n \geq 5$. By Lemma \ref{l1},
$S = \{f(u) + f(v) : uv \in E(H_n)\}$ is a set of  $|E(H_n)| = 2|V(H_n)|
- 3$ consecutive integers. Thus $S  = \{3, 4, 5, \ldots, 2n, 2n+1\}$.
We shall consider two cases.

\textbf{Case 1}: $n = 5, 6$. For $n = 5$, The sum of all elements of $S$ is 63. This sum contains
two times of label $x_1$, three times each label of $x_i$, $2 \leq
i \leq 5$  and four times the label of $c$. Thus, we have
\[\sum_{i=2}^5 f(x_i)+ 2f(c) = 21.\]
It is a routine procedure to verify that this equation has no
solution. Hence, $H_5$ is not super edge-magic. With a similar
argument, for $n = 6$, we have
\[\sum_{i=2}^6 f(x_i)+ 3f(c) = 32.\]
The possible solutions for this equation are $f(c) = 3, f(x_1) = 2$,  $f(x_i) \in \{1, 4, 5, 6, 7\}, 2 \leq i \leq 6,$ and
$f(c) = 5, f(x_1) = 6$,  $f(x_i) \in \{1, 2, 3,  4, 7\}, 2 \leq i \leq 6$.  It can be checked that these solutions do not lead to a super edge-magic labeling of $H_6$. Hence, $H_6$ is not a super edge-magic graph.

\textbf{Case 2}: $n \geq 7$. Observe that both 3 and 4 can be expressed uniquely as sums of
two distinct element from the set $\{1, 2, 3, \ldots, n+1\}$, namely
$3= 1+2$ and $4 = 1+3$. On the other hand, 5 can be expressed as sums of distinct elements of $\{1, 2, \ldots,
n+1\}$ in exactly two ways, namely $5 = 2 + 3 = 1 + 4$. Then, the vertices of labels 1,
2 and 3 must form a triangle or the vertex of label 1 is adjacent to
the vertices of labels 2, 3 and 4, respectively.
%Also, the only possibility to get the sums $2n+1$, $2n$ and $2n-1$ in S is
%$2n+1 = (n+1) + n$, $2n = (n+1) + (n-1)$ and $2n-1 = n + (n-1) = (n+1) + (n-2)$.
With a similar argument, the vertices of labels $n-1$, $n$ and $n+1$
must form a triangle or the vertex of label $n+1$ is adjacent to the
vertices of labels $n$, $n-1$ and $n-2$, respectively. By combining these facts, we obtain either $2K_3$, $K_3 \cup K_{1,3}$ or
$2K_{1,3}$ as a subgraph of $H_n$, a contradiction. This
completes the proof. \BOX

Based on the results of Theorem \ref{T3}, the super edge-magic deficiency of
$H_n$ is 0 for $n = 3$ and 4, and at least 1 for $n \geq 5$. For
$n = 5, 6, 7$, we could prove that $\mu_s(H_n) = 1$ by labeling the vertices $(c; x_1, x_2, x_3, x_4, x_5)$, $(c; x_1, x_2,
x_3, x_4, x_5, x_6)$, and $(c; x_1, x_2, x_3, x_4, x_5, x_6, x_7)$
with $(1; 7, 5, 3, 6, 4)$, $(2; 3, 1, 4, 8, 5, 6)$, and $(2; 3, 1,
4, 8, 5, 9, 6)$, respectively.

For $n \geq 8$ we shall determine an upper bound for the super edge-magic deficiency of $H_n$ where $n \neq 2 \mod 4$ as stated in the following
theorem.

\begin{theorem}
For any integer $n \geq 8$, $n \equiv 0, 1, 3$ (mod 4), the super edge-magic deficiency of $H_n$ are given
by
\[\mu_s(H_n) \leq \left\{\begin{array}{ll} \frac{1}{2}(n-3), &\mbox{if}~n \equiv 1~\mbox{or}~3~\mbox{mod~4},\\
\frac{n}{2}, &\mbox{if}~n \equiv 0~\mbox{mod~4}.
\end{array}\right.\]
\end{theorem}
\pr We consider the following two cases.

\textbf{Case 1}: $n \equiv 1$ or 3 (mod 4). Define a vertex labeling as follow.
\[f(c) = \frac{1}{2}(3n-1).\]
\[f(x_i) = \left\{\begin{array}{ll}  \frac{1}{2}(i+1),  &\mbox{if}~i = 1, 3, 5, \ldots, n-1,\\
\lceil\frac{n}{2}\rceil + \frac{i}{2}, &\mbox{if}~i = 2, 4, 6,
\ldots, n-2.
\end{array}\right.\]

\textbf{Case 2}: $n \equiv 0$  (mod 4). We redefine the edge-set of $H_n$ as $E(H_n) = \{cx_i : 1 \leq i \leq n, i \neq \frac{n}{2}\} \cup \{x_ix_{i+1} : 1 \leq i \leq n-1\} \cup \{x_nx_1\}$.
Now we are ready to define a vertex labeling $f$.
\[f(c) = \frac{1}{2}(3n+2).\]
\[f(x_i) = \left\{\begin{array}{ll}  \frac{1}{2}(i+1),  &\mbox{if}~i = 1, 3, 5, \ldots, n-1,\\
 \frac{1}{2}(n+i), &\mbox{if}~i= 2, 4, 6, \ldots,\frac{1}{2}(n-4),\\
\frac{5}{4}n, & \mbox{if}~i= \frac{n}{2},\\
\frac{1}{2}(n+i-2), &\mbox{if}~i = \frac{n}{2}+2, \frac{n}{2}+4,
\ldots, n.
\end{array}\right.\]
For both cases, it is easy to verify that $f$ extends to a super edge-magic labeling of $H_n$. \BOX

We have tried to find an upper bound of the super edge-magic
deficiency of $H_n$ for $n \equiv 2$ (mod 4), but without
success. And thus we propose the following problems.

\begin{open problem}
For $n \equiv 2$ (mod 4), find an upper bound of the super
edge-magic deficiency of $H_n$. Further, find the super
edge-magic deficiency of $H_n$ for all $n$.
\end{open problem}
%---------------------------------------------------------

\section{Super edge-magic deficiency of join-product graphs}

In this section, we consider super edge-magic deficiency of three classes of graphs.
These graphs are obtained from join products of a path $P_n$, a star $K_{1,n}$, and a cycle $C_n$, respectively, with $m$ isolated vertices ($\overline{K_m}$).

First, we consider the super edge-magic deficiency of $P_n + \overline{K_m}$. We denote the vertex and edge sets of $P_n + \overline{K_m}$ as
\[V(P_n + \overline{K_m}) = \{u_i : 1 \leq i \leq n\} \cup \{v_j : 1 \leq j \leq
m\}\] and
\[E(P_n + \overline{K_m}) = \{u_iu_{i+1} : 1 \leq i \leq n-1\} \cup \{u_iv_j : 1 \leq i \leq
n, 1 \leq j \leq m\}.\] It is clear that $P_n + \overline{K_m}$ is a graph of
order and size $n+m$ and $n(m + 1)-1$, respectively.

If $m = 1$, then $P_n + \overline{K_1}$ is a fan $F_n$. As we
mention in the first section, the super edge-magic deficiency of $F_n$ have been
studied in \cite{N1}. Furthermore, Ngurah \textit{et al.} \cite{N2} studied the super edge-magic
deficiency of $P_n + \overline{K_2}$ and proved that $\mu_s(P_n + \overline{K_2}) =
\frac{1}{2}(n-2)$ for all even $n \geq 2$, and conjectured that
$\mu_s(P_n + \overline{K_2}) = \frac{1}{2}(n-1)$  for all odd $n \geq 3$. In this
section, we study the super edge-magic deficiency of $P_n + \overline{K_m}$
for $m \geq 3$. The next result provides sufficient and necessary
conditions for $P_n + \overline{K_m}$ to be super edge-magic.

\begin{lemma}\label{Lm1}
Let $n \geq 1$ and $m \geq 3$ be integers. Then the graph $P_n + \overline{K_m}$ is super edge-magic if and only if $n \in \{1, 2\}.$
\end{lemma}
\pr First, we show that $P_n + \overline{K_m}$ is super edge-magic for $n = 1,
2$. It is known that $P_1 + \overline{K_m} \cong K_{1, m}$ is super edge-magic. For $n = 2$ label the vertices
$\{u_1, u_2\}$ and $\{v_1, v_2, v_3, \ldots, v_m\}$ with $\{1, m+2\}$ and $\{2, 3, \ldots, m+1\}$, respectively. Then by Lemma
\ref{l1}, this vertex labeling extends to a super edge-magic labeling of $P_2 + \overline{K_m}$ with the magic constant $3m+6$. For the sufficiency, let $P_n + \overline{K_m}$ be a super edge-magic graph. By Lemma \ref{l2}, we have $n(m + 1)-1 \leq 2(n+m) - 3$ and the desired result.\BOX

Based on Lemma \ref{Lm1},  $\mu_s(P_n + \overline{K_m}) = 0$ for $n
\leq 2$ and $\mu_s(P_n + \overline{K_m}) \geq 1$ for $n \geq 3$. Since there is
no super edge-magic labeling of $P_n + \overline{K_m}$ for almost all values of
$n$, we thus try to find its super edge-magic deficiency. The
following theorem gives upper and lower bounds of the
deficiency.

\begin{theorem} \label{T1}
For any integers $n, m \geq 3$, the super edge-magic
deficiency of $P_n + \overline{K_m}$ satisfies $\lceil
\frac{1}{2}(n-2)(m-1)\rceil \leq \mu_s(P_n + \overline{K_m}) \leq  (n-1)(m-1)-1.$
\end{theorem}
\pr To prove the upper bound, we define a vertex labeling $f$ as follow.
\[f(u_i)=\left\{\begin{array}{ll} \lfloor\frac{1}{2}(n + 2)\rfloor + \frac{1}{2}(i-1), &\mbox{for~ood}~ i,\\
n + \frac{1}{2}i, &\mbox{for~even}~ i,\\
\end{array}\right.\]
and
\[f(\{v_1, v_2, v_3, \ldots, v_m\}) = \{1, 2n, 3n, 4n, \ldots,
mn\}.\]
We can see that these vertex-labels are non-repeated and constitute a set
$\{f(x) + f(y) | xy \in E(P_n + \overline{K_m})\}$ of $n(m+1)-1$ consecutive
integers. However, the largest vertex label used is $mn$ and there exist $mn-(n+m) =
(n-1)(m-1)-1$ labels that are not utilized. So, for each number between 1 and $mn$ that has not been used as a label, we
introduce a new vertex labeled with that number; and this gives
$(n-1)(m-1)-1$ isolated vertices. By Lemma \ref{l1}, this yields
a super edge-magic  labeling of  $P_n + \overline{K_m} \cup [(n-1)(m-1)-1] K_1$
with magic constant $2mn + \lfloor\frac{1}{2}(3n + 2)\rfloor$.
Hence,
\[\mu_s(P_n + \overline{K_m}) \leq (n-1)(m-1)-1.\]
For a lower bound, by Lemma \ref{l2}, it is easy to verify that
\begin{center}
$\mu_s(P_n + \overline{K_m}) \geq \lceil \frac{1}{2}(n-2)(m-1)\rceil$. \BOX
\end{center}

%For a lower bound, by Lemma \ref{l2}, it is easy to check that
%$F_{n, m} \cup [\lceil \frac{1}{2}(n-2)(m-1)\rceil - 1] K_1$ is not
%super edge-magic. Thus,
%\begin{center}
%$\mu_s(F_{n,2}) \geq \lceil \frac{1}{2}(n-2)(m-1)\rceil$. \BOX
%\end{center}
Notice that, the lower bound presented in Theorem \ref{T1} is
tight. We  found that the super edge-magic deficiency of  $P_4 + \overline{K_m}$ is equal to its lower bound by
labeling the vertices $(u_1,  u_2, u_3, u_4)$ and $\{v_1, v_2,  v_3,
\ldots, v_m\}$ with $(1, 2, 2m+2, 2m+3)$ and $\{3, 5, 7, \ldots,
2m-1, 2m+1\}$, respectively. This vertex-labels extend to a super
edge-magic labeling of $P_4 + \overline{K_m}$ with the magic constant $6m + 9$.
The largest vertex label used is $2m+3$. So, $\mu_s(P_4 + \overline{K_m}) \leq
2m+ 3-(m+4) = m-1.$ From this fact and Theorem \ref{T1},
$\mu_s(P_4 + \overline{K_m}) = m-1$. Additionally, we found that $\mu_s(P_6 + \overline{K_m}) =
2(m-1)$ by labeling the vertices $(u_1, u_2, u_3, u_4, u_5, u_6)$ and
$\{v_1, v_2, v_3, \ldots, v_m\}$ with $(2, 1, 3, 3m+2, 3m+4, 3m+3)$
and $\{4, 7, 10, \ldots, 2m-5, 2m-2, 3m+1\}$, respectively.
%This
%vertex-labels extend to a super edge-magic labeling of $P_6 \odot \overline{K_m}$
%with the magic constant $9m + 12$.

Referring to the afore-mentioned results, we propose the following problems.
\begin{open problem}
Find a better upper  bound of the super edge-magic deficiency of $P_n + \overline{K_m}$. Further, find the  super
edge-magic deficiency of $P_n + \overline{K_m}$ for $n \neq 4, 6$.
\end{open problem}

Let us now determine the super edge-magic deficiency of  $K_{1,n} + \overline{K_m}$. Let $K_{1,n} + \overline{K_m}$ be a graph having
\[V(K_{1,n} + \overline{K_m}) = \{c\} \cup \{x_i : 1 \leq i \leq n\} \cup \{y_j : 1 \leq j \leq m\},\]
and
\[E(K_{1,n} + \overline{K_m}) = \{cx_i : 1 \leq i \leq n\} \cup \{x_iy_j, cy_j : 1 \leq i \leq n, 1 \leq j \leq m\}.\]
Thus $K_{1,n} + \overline{K_m}$ is a graph of order $n+m+1$ and size $(n+1)(m+1)-1.$ Notice
that if $n = 1$, then $K_{1,1} + \overline{K_m} \cong P_2 + \overline{K_m}$ which
is super edge-magic (see Theorem \ref{T1}). Hence, we assume that $n \geq 2.$

\begin{lemma}\label{L2} Let $n \geq 2$ and $m \geq 1$ be integers. Then,
$K_{1,n} + \overline{K_m}$ is super edge-magic if and only if $m =
1.$
\end{lemma}
\pr By Lemma \ref{l2}, it is easy to check that if $K_{1,n} + \overline{K_m}$ is super
edge-magic then $m \leq 1.$ Since $m$ is a positive integer, so $m = 1.$  For the sufficiency, label the vertices $\{c\},$
$\{x_1, x_2, x_3, \ldots, x_n\},$ and $\{y_1\}$ with $\{n+1\},$
$\{1, 2, 3, \ldots, n\},$ and $\{n+2\}$, respectively. This vertex
labeling extends to a super edge-magic labeling of $K_{1,n} + \overline{K_m}$ with
magic constant $3n+6$. \BOX

Since $K_{1,n} + \overline{K_m}$ is not super edge-magic for almost all values of $m$, we
thus try to find its super edge-magic deficiency. The following
result gives upper and lower bounds of the deficiency.

\begin{theorem} \label{T2}
For any integers $n, m \geq 2$, the super edge-magic
deficiency of $K_{1,n} + \overline{K_m}$ satisfies $\lceil
\frac{1}{2}(n-1)(m-1)\rceil \leq \mu_s(K_{1,n} + \overline{K_m}) \leq n(m-1)-1.$
\end{theorem}
\pr Similar with the proof of Theorem \ref{T1}, we could obtain
that $\mu_s(K_{1,n} + \overline{K_m}) \geq \lceil \frac{1}{2}(n-1)(m-1)\rceil.$ To show the
upper bound, label the vertices $\{c\},$ $\{x_1, x_2, x_3, \ldots,
x_n\},$ and $\{y_1, y_2, y_3, \ldots, y_m\}$ with $\{n+2\},$ $\{2,
3, 4, \ldots, n+1\},$ and $\{1, 2(n+1), 3(n+1), \ldots, m(n+1)\}$,
respectively. This vertex labeling extends to a super edge-magic
labeling of $K_{1,n} + \overline{K_m}$ with magic constant $(n+1)(m+1)+1$ and the
largest vertex label $m(n+1).$ \BOX

\begin{open problem}
For integers $n, m \geq 2$, find better upper and lower bounds of the super
edge-magic deficiency of $K_{1,n} + \overline{K_m}$. Further, find the super
edge-magic deficiency of $K_{1,n} + \overline{K_m}$ for a fixed value of $n$ or $m$.
\end{open problem}

Finally, we consider the super edge-magic deficiency of $C_n + \overline{K_m}$. Notice that this graph is not super edge-magic for all integers $n \geq 3$ and $m \geq 1$. For $m = 1$, the graph $C_n + \overline{K_1}$ is a wheel $W_{n}$.
Ngurah \textit{et al.} \cite{N1} studied the super edge-magic deficiency of $W_{n}$ and they determined the super edge-magic deficiency of $W_{n}$
for some values of $n$ and gave a lower bound for general values of $n$. Additionally, they also provided an upper bound for the super edge-magic deficiency of $W_{n}$ for odd $n$. Now, we study the super edge-magic deficiency of $C_n + \overline{K_m}$ for $n \geq 3$ and $m \geq 2$. Our first result gives a lower bound of the super edge-magic deficiency of $C_n + \overline{K_m}$.

\begin{lemma}\label{L3} For any integers $n \geq 3$ and $m \geq 2$,
$\mu_s(C_n + \overline{K_m}) \geq \lfloor\frac{1}{2}(m+1)n\rfloor - (n+m) + 2.$
\end{lemma}
\pr It is easy to verify that $C_n + \overline{K_m} \cup tK_1$, where $t = \lfloor\frac{1}{2}(m+1)n\rfloor - (n+m) + 1$, is not a super edge-magic graph. Hence, $\mu_s(C_n + \overline{K_m}) \geq \lfloor\frac{1}{2}(m+1)n\rfloor - (n+m) + 2.$ \BOX

\begin{theorem} \label{T4}
Let $n \geq 3$ be an odd integer. Then $\mu_s(C_n + \overline{K_m}) \leq mn - (n+m) + 1$ for every  integer $m \geq 2$.
\end{theorem}
\pr Let $C_n + \overline{K_m}$ be a graph with
$$V(C_n + \overline{K_m}) = \{u_i : 1 \leq i \leq n\} \cup \{v_i : 1 \leq j \leq m\}$$
and
$$E(C_n + \overline{K_m}) = \{u_iu_{i+1} : 1 \leq i \leq n-1\} \cup \{u_nu_1\} \cup \{u_iv_j : 1 \leq i \leq n, 1 \leq j \leq m\}.$$

Next, define a vertex labeling $f$ as follow.
\[f(u_i) = \left\{\begin{array}{ll}  \frac{1}{2}(n+2+i),  &\mbox{if}~i = 1, 3, 5, \ldots, n,\\
 \frac{i}{2}(2n+2+i), &\mbox{if}~i = 2, 4, 6,
\ldots, n-1,
\end{array}\right.\]
$$f(\{v_1, v_2, v_3, \ldots, v_m\}) = \{1, 2n+1, 3n+1, 4n+1, \ldots, mn+1\}.$$

It is a routine procedure to check that $f$  can be extended to a super edge-magic labeling of
$C_n + \overline{K_m} \cup tK_1$, where $t = mn - (n+m) + 1.$ Thus, we have the desired result.  \BOX

Some open problems related the super edge-magic deficiency of $C_n + \overline{K_m}$ are
presented bellow.
\begin{open problem}
For even $n \geq 4$ and every $m \geq 2$, find an upper bound for the super
edge-magic deficiency of $C_n + \overline{K_m}$. Further, find a better upper bound of the super edge-magic deficiency of $C_n + \overline{K_m}$ for odd $n$ and
every $m \geq 2$.
\end{open problem}

%So far, we studied the super edge-magic deficiency of  join product of particular type of super edge-magic graphs (a path, a star, and a cycle, respectively), with isolated vertices. As we can see,  these  graphs have finite super edge-magic deficiency. In fact, the join products of any  super edge-magic graph with isolated vertices has finite super edge-magic deficiency, as we state in the following theorem.
Our results showed the finiteness of super edge-magic deficiencies of join product of a path, a star, and a cycle with isolated vertices.  Recall that all paths, stars, and cycles of odd order are super edge-magic.
In the next theorem, we managed to generalize similar result for any super edge-magic graph.
\begin{theorem}
Let  $G$ be a super edge-magic graph with a super edge-magic labeling $f$. For any integer $m \geq 1,$  $\mu_s(G + \overline{K_m}) \leq s + (m-2)|V(G)|-m,$ where $s = max\{f(u) + f(v) : uv \in E(G)\}.$
\end{theorem}

%\begin{theorem}
%If $G$ is a super edge-magic graph, then  the super edge-magic
%deficiency of $G + \overline{K_m}$ is finite for any integer $m \geq 1$.
%\end{theorem}
\pr First, define $H \cong G + \overline{K_m}$ as a graph with $V(H) = V(G) \cup \{y_1, y_2, y_3, \ldots, y_m\}$ and $E(H) = E(G) \cup \{xy_i : x \in V(G), 1 \leq i \leq m\}.$
Next, define a vertex labeling $g$ as follows. $$g(x) = f(x),~\mbox{if}~x \in V(G),$$ and $$g(\{y_1, y_2, y_3, \ldots, y_m\}) = \{s, s+|V(G)|, s+2|V(G)|, \ldots, s+(m-1)|V(G)|\}.$$
It is easy to verify that $g$ extends to a super edge-magic labeling of
$H \cup [s + (m-2)|V(G)|-m]K_1$. Hence, $\mu_s(H) \leq s + (m-2)|V(G)|-m.$ \BOX

To conclude, we would like to ask an interesting general question regarding the super edge-magic deficiency of join-product graphs.
\begin{open problem}
If $G$ is an arbitrary graph, determine the super edge-magic deficiency of the join-product of $G$ with $m$ isolated vertices, $\mu_s(G + \overline{K_m})$.
\end{open problem}

\end{document}